\documentclass[reqno]{amsart}

\usepackage{amsmath, amsfonts, amssymb, hyperref}
\usepackage{graphics}
\usepackage[usenames]{color}

{\theoremstyle{plain}
\newtheorem{theorem}{Theorem}[section]

\newtheorem{proposition}[theorem]{Proposition}
\newtheorem{lemma}[theorem]{Lemma}
\newtheorem{corollary}[theorem]{Corollary}
}
{\theoremstyle{definition}
\newtheorem{definition}[theorem]{Definition}
\newtheorem*{definition*}{Definition}
\newtheorem{remark}[theorem]{Remark}
\newtheorem{example}[theorem]{Example}

}

\numberwithin{equation}{section}

\newcommand{\domleq}{\preceq}   
\newcommand{\supp}[1]{\mathrm{supp}(#1)}
\newcommand{\rows}[2]{\mathrm{rows}_{#1}(#2)}
\newcommand{\cols}[2]{\mathrm{cols}_{#1}(#2)}
\newcommand{\rects}[3]{\mathrm{rects}_{#1,#2}(#3)}
\newcommand{\trim}[2]{\mathrm{trim}^{#1}(#2)}
\newcommand{\union}{\cup}
\newcommand{\overlap}[2]{\mathrm{overlap}_{#1}(#2)}

\begin{document}
\title[Necessary Conditions for Schur-Positivity]{Necessary Conditions for Schur-Positivity}

\author{Peter R. W. McNamara}
\address{Department of Mathematics, Bucknell University, Lewisburg, PA 17837, USA}
\email{\href{mailto:peter.mcnamara@bucknell.edu}{peter.mcnamara@bucknell.edu}}


\begin{abstract} 
In recent years, there has been considerable interest in
showing that certain conditions on skew shapes $A$ and $B$ are sufficient for
the difference $s_A - s_B$ of their skew Schur functions
to be Schur-positive.  
We determine \emph{necessary} conditions
for the difference to be Schur-positive.  Our conditions are motivated by those 
of Reiner, Shaw and van Willigenburg that are necessary for $s_A = s_B$, and we deduce a
strengthening of their result as a special case.  

\end{abstract}

\maketitle


\section{Introduction}\label{sec:intro} 

In many respects, the basis of Schur functions is the most interesting and important basis for the
ring of symmetric functions.
The significance of Schur functions is highlighted by their appearance in several areas of
mathematics.  In particular, they arise in the representation theory of the symmetric group
and of the general and special linear groups.  They appear in algebraic geometry, specifically
in the study of the cohomology ring of the Grassmannian, and they are also closely connected
to the eigenvalues of Hermitian matrices. 

It is therefore natural to study the expansion of symmetric functions as a linear combination
of Schur functions.   For example, skew Schur functions $s_{\lambda/\mu}$ and the product 
$s_\sigma s_\tau$ of two 
Schur functions are famous examples of \emph{Schur-positive} functions:
they can be written as linear combinations of Schur functions with all coefficients positive.
Taking this a step further, several recent papers such as 
\cite{BBR06, FFLP05, Kir04, LPP07, LLT97, Oko97}
have asked when expressions of the form 
\[
s_{\lambda/\mu} - s_{\sigma/\tau} \mbox{\ \ \ \ or\ \ \ \ } s_\lambda s_\mu - s_\sigma s_\tau
\]
are Schur-positive.  These papers have been concerned with giving conditions
on $\lambda, \mu, \sigma$ and $\tau$ that result in Schur-positive expressions.  
We wish to focus on the converse direction:  if we know that the expressions are Schur-positive, 
what must be true about $\lambda, \mu, \sigma$ and $\tau$?


Let us first note that $s_\lambda s_\mu$ is just a special type of skew Schur function 
(see Subsection~\ref{sub:products} for an explanation).  Therefore, it suffices to 
consider differences of the form $s_A - s_B$, where $A$ and $B$ are skew shapes.
It is well-known that if $s_A - s_B$ is Schur-positive, then the partition of row lengths
of $B$ must dominate the partition of row lengths of $A$; 
see Proposition~\ref{pro:extreme_fillings}.  Similarly, the partition of column lengths
of $B$ must dominate the partition of column lengths of $A$.  
In \cite{RSV07}, some necessary conditions on $A$ and $B$ are given for the equality
$s_A = s_B$ to hold; these conditions depend not only on the rows lengths of $A$ and
$B$, but also on the overlaps between the various rows.  Inspired by this, our main result, 
Corollary~\ref{cor:combine}, roughly
says that if $s_A - s_B$ is Schur-positive, then all the row overlaps for $B$ must dominate
those of $A$.  It is worth mentioning that this result includes the well-known results
about the partitions of row lengths and column lengths as special cases.
The full details are the content of Section~\ref{sec:main}.

In fact, our results require a weaker condition than the Schur-positivity of $s_A - s_B$.  We say that the support
of a skew shape $A$ is the set of partitions $\lambda$ such that $s_\lambda$ appears with nonzero
coefficient when we expand $s_A$ in terms of Schur functions.  
Instead of requiring that $s_A - s_B$ is Schur-positive, our proofs only require
that the support of $A$ contains the support of $B$.  This allows us to strengthen the
aforementioned result
of \cite{RSV07}, which we do in the first part of Section~\ref{sec:specialcases}.
In the rest of Section~\ref{sec:specialcases}, we restrict our results to obtain necessary
conditions for $s_\lambda s_\mu - s_\sigma s_\tau$ to be Schur-positive.  

\subsection*{Acknowledgements}  The author
 thanks Stephanie van Willigenburg for carefully reading and
giving useful comments on an early version of this manuscript.  Both \cite{BucSoftware} and
\cite{SteSoftware} aided invaluably in data generation.

\section{Preliminaries}\label{sec:prelims} 

We follow the terminology and notation of \cite{Mac95} and \cite{ECII}.

\subsection{Skew shapes} 
A \emph{partition} $\lambda$ of $n$ is a weakly decreasing list of positive integers 
$(\lambda_1, \ldots, \lambda_l)$ whose sum is $n$.  We say that $n$ is the size of $\lambda$, 
denoted $|\lambda|$, and we call $l$ the \emph{length} of
$\lambda$ and denote it by $\ell(\lambda)$.  It will be convenient to set $\lambda_k = 0$
for $k > \ell(\lambda)$, thus identifying $\lambda$ with
$(\lambda_1, \ldots, \lambda_{\ell(\lambda)}, 0, 0, \ldots, 0)$, where the string of zeros has
arbitrary length.  In particular, the unique partition of $0$ can be denoted by $(0)$.  
We will mainly think of $\lambda$ in terms of its \emph{Young diagram}, which is a
left-justified array of boxes that has $\lambda_i$ boxes in the $i$th row from the top.  For
example, if $\lambda = (4,4,3)$, which we will abbreviate as $\lambda = 443$, then the 
Young diagram of $\lambda$ is 
\setlength{\unitlength}{4mm}
\[
\begin{picture}(4,3)(0,0)
\put(0,3){\line(1,0){4}}
\put(0,2){\line(1,0){4}}
\put(0,1){\line(1,0){4}}
\put(0,0){\line(1,0){3}}
\multiput(0,3)(1,0){5}{\line(0,-1){2}}
\multiput(0,1)(1,0){4}{\line(0,-1){1}}
\end{picture}.
\]
We will say that a partition $\mu$ is \emph{contained} in a partition $\lambda$ if the
Young diagram of $\mu$ is contained in the Young diagram of $\lambda$.  
In this case, we define the \emph{skew shape} $\lambda/\mu$ to be the set of boxes in
the Young diagram of $\lambda$ that remain after we remove those boxes corresponding
to $\mu$.  For example, the skew shape $A = (4,4,3)/(2)$ is represented as
\[
\begin{picture}(4,3)(0,0)
\put(2,3){\line(1,0){2}}
\put(0,2){\line(1,0){4}}
\put(0,1){\line(1,0){4}}
\put(0,0){\line(1,0){3}}
\multiput(2,3)(1,0){3}{\line(0,-1){1}}
\multiput(0,2)(1,0){5}{\line(0,-1){1}}
\multiput(0,1)(1,0){4}{\line(0,-1){1}}
\end{picture}.
\]
We will label skew shapes by simply
using single uppercase roman letters, as in the example above.  We 
write $|A|$ for the \emph{size} of $A$, which is simply the number of boxes in the skew shape $A$.
If $A = \lambda/\mu$ and $\mu=(0)$, then $A$ is said to be a \emph{straight shape}.

\subsection{Skew Schur functions and the Littlewood-Richardson rule}
While skew shapes are our main diagrammatical objects of study, our main algebraic
objects of interest are skew Schur functions, which we now define.  For a skew shape $A$, a
\emph{semi-standard Young tableau} (SSYT) of shape $A$ is a filling of the boxes of
$A$ with positive integers such that the entries weakly increase along the rows and 
strictly increase down the columns.  For example, 
\[
\begin{picture}(4,3)(0,0)
\put(2,3){\line(1,0){2}}
\put(0,2){\line(1,0){4}}
\put(0,1){\line(1,0){4}}
\put(0,0){\line(1,0){3}}
\multiput(2,3)(1,0){3}{\line(0,-1){1}}
\multiput(0,2)(1,0){5}{\line(0,-1){1}}
\multiput(0,1)(1,0){4}{\line(0,-1){1}}
\put(2.3,2.2){1}
\put(3.3,2.2){2}
\put(0.3,1.2){1}
\put(1.3,1.2){1}
\put(2.3,1.2){2}
\put(3.3,1.2){3}
\put(0.3,0.2){5}
\put(1.3,0.2){7}
\put(2.3,0.2){7}
\end{picture}
\]
is an SSYT of shape $443/2$. 
The \emph{skew Schur function} $s_A$ in the variables $(x_1, x_2, \ldots)$ is then defined by
\[
s_A = \sum_{T} x^T
\]
where the sum is over all SSYT $T$ of shape $A$, and  
\[
x^T = x_1^{\mbox{\scriptsize \#1's in $T$}}
x_2^{\mbox{\scriptsize \#2's in $T$}} \cdots .
\]
For example, the SSYT above contributes the monomial $x_1^3 x_2^2 x_3 x_5 x_7^2$ to 
$s_{443/2}$.  The sequence $(\mbox{\#1's in $T$, \#2's in $T$, $\ldots$})$ is known as
the \emph{content} of $T$.

If $A$ is a straight shape, then $s_A$ is called simply a \emph{Schur function},
and some of the significance of Schur functions stems from the fact that they form a basis for the
symmetric functions.  Therefore, every skew Schur function can be written as a linear
combination of Schur functions.  A simple description of the coefficients in 
this linear combination is given by the celebrated \emph{Littlewood-Richardson rule}, which we 
now describe.  The \emph{reverse reading word} of an SSYT $T$ is the word obtained by 
reading the entries of $T$ from right to left along the rows, taking the 
rows from top to bottom.  For example, the SSYT above has reverse reading word 213211775.
An SSYT $T$ is said to be an \emph{LR-filling} if, as we read the reverse reading word of $T$,
the number of appearances of $i$ always stays ahead of the number of appearances of
$i+1$, for $i=1,2,\ldots$.
The reader is invited to check that the only possible LR-fillings of $443/2$ have
reading words 112211322 and 112211332.  The Littlewood-Richardson rule \cite{LiRi34, Sch77, ThoThesis, Tho78} then states that
\[
s_{\lambda/\mu} = \sum_{\nu} c^\lambda_{\mu\nu} s_\nu\ ,
\]
where $c^\lambda_{\mu\nu}$ is the ubiquitous \emph{Littlewood-Richardson} coefficient, 
defined to be the number of LR-fillings of $\lambda/\mu$ with content $\nu$. 
For example, if $A = 443/22$, then $s_A = s_{441} + s_{432}$.
It follows that any skew Schur function can be written 
as a linear combination of Schur functions with all positive coefficients, and we thus say
that skew Schur functions are \emph{Schur-positive}.  

As mentioned in the introduction,
our main goal is to determine when the difference $s_A - s_B$ of two skew Schur functions
is Schur-positive.  It turns out that most of our results can be expressed in terms of the 
\emph{support} of skew Schur functions.  The support $\supp{A}$ of $s_A$ 
is defined to be the set of those partitions $\nu$ for
which $s_\nu$ appears with nonzero coefficient when we expand $s_A$ in terms of 
Schur functions.  For example, we have $\supp{443/2} = \{441, 432\}$.  

We will make significant use of the transpose operation on skew shapes and we will
also need the related $\omega$ involution on symmetric functions.  For any
partition $\lambda$, we define the transpose $\lambda^t$ to be the partition obtained by reading
the column lengths of $\lambda$ from left to right.  For example, $(443)^t = 3332$.
The transpose operation can be extended to skew shapes $A=\lambda/\mu$ by setting
$A^t = \lambda^t/\mu^t$.  Then $\omega$ is defined by $\omega(s_\lambda) = s_{\lambda^t}$.  
It can be shown that $\omega(s_A) = s_{A^t}$ for any skew shape $A$.

\subsection{Extended dominance order} 
The well-known dominance order is typically restricted to partitions of equal size, but
its definition readily extends to give
a partial order on arbitrary partitions, and this is our final
preliminary.

\begin{definition}
For partitions $\lambda = (\lambda_1, \lambda_2, \ldots, \lambda_r)$ and 
$\mu = (\mu_1, \mu_2, \ldots, \mu_s)$, we define the \emph{dominance order} $\domleq$ by
$\lambda \domleq \mu$ if 
\[
\lambda_1 + \lambda_2 + \cdots \lambda_k \leq \mu_1 + \mu_2 + \cdots \mu_k
\]
for all $k=1,2,\ldots,r$. 
In this case, we will say that $\mu$ \emph{dominates} $\lambda$, or is \emph{more dominant}
than $\lambda$.
\end{definition}

Note that our definition does not require that $|\lambda| = |\mu|$.  For example, we have
$(4,2,1) \domleq (4,4)$.

We will need the following result about our extended definition of dominance order.  Since
it is straightforward to check, we leave the proof as an exercise.

\begin{lemma}\label{lem:dom_containment}
Consider two sequences $a = (a_1, a_2, \ldots, a_r)$ and $b = (b_1, b_2, \ldots, b_s)$ of natural
numbers such that 
$r \leq s$ and $a_i \leq b_i$ for $i=1,2,\ldots,r$.  Let $\alpha$ and $\beta$ denote the partitions
obtained by sorting the parts of $a$ and $b$ respectively into weakly decreasing order.
Then $\alpha \domleq \beta$.  
\end{lemma}

\section{Necessary conditions for Schur-positivity}\label{sec:main}

We begin in earnest by stating well-known necessary conditions for $s_A - s_B$ to be Schur-positive,
which nevertheless seem to be absent from the literature.
Since necessary conditions are our
focus and we wish to make our presentation self-contained, we will also give a proof.
For any skew shape $A$, let $\rows{}{A}$ denote the partition obtained by sorting the 
row lengths of $A$ into weakly decreasing order.  Similarly, let $\cols{}{A}$ be the partition
obtained from the column lengths.  

\begin{proposition}\label{pro:extreme_fillings}  Let $A$ and $B$ be skew shapes.
If $\lambda \in \supp{A}$, then 
\[
\rows{}{A} \domleq \lambda \domleq \cols{}{A}^t
\]
and both inequalities are sharp.
Consequently, if $s_A - s_B$ is Schur-positive, then
\[
\rows{}{A} \domleq \rows{}{B}  \mbox{\ \ and \ \ } \cols{}{A} \domleq \cols{}{B}.
\]
\end{proposition}

\begin{figure}[htpb]
\[
\scalebox{.75}{\includegraphics{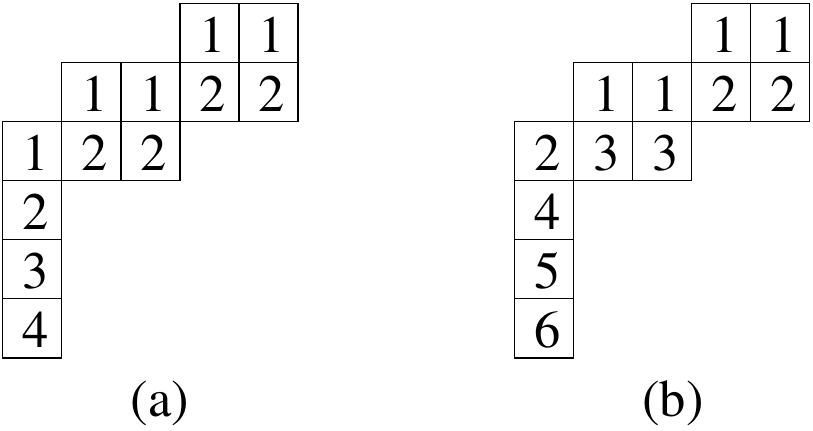}}
\]
\caption{The most and least dominant fillings of $553111/31$.}
\label{fig:extreme_fillings}
\end{figure}

\begin{proof}
We first show that $\lambda \domleq \cols{}{A}^t$.  Suppose we wish to construct an LR-filling
of $A$ that is as dominant as possible.  Thus we wish to use as many 1's as possible, then 
as many 2's as possible, and so on.  Since we can only have at most one $i$ in each column of
$A$, fill the $i$th highest box of every column of $A$ with the number $i$, for all $i$.  It is 
straightforward to check that the result is an LR-filling of $A$ of content $\cols{}{A}^t$, and 
that every other LR-filling of $A$ will have a less dominant content.  
See Figure~\ref{fig:extreme_fillings}(a) for an example.  

Applying the $\omega$ involution, we know that $\lambda \in \supp{A}$ if and only if 
$\lambda^t \in \supp{A^t}$.  Thus $\lambda^t \domleq \cols{}{A^t}^t = \rows{}{A}^t$.  
As shown in \cite{Bry73}, when partitions $\mu$ and $\nu$ satisfy $|\mu|=|\nu|$, the transpose 
operation is order-reversing with respect to $\domleq$\:; i.e. $\mu \domleq \nu$ if and 
only if $\nu^t \domleq \mu^t$.  We conclude that $\rows{}{A} \domleq \lambda$.
This inequality is sharp since $\lambda^t \domleq \cols{}{A^t}^t $ is sharp.

If $s_A - s_B$ is Schur-positive, then $\supp{s_B} \subseteq \supp{s_A}$ and so 
$\rows{}{A} \domleq \rows{}{B}$ and $\cols{}{B}^t \domleq \cols{}{A}^t$.  The result now
follows from the order-reversing property of the transpose operation.
\end{proof}

\begin{remark}\label{rem:extreme_fillings}
It follows from Proposition~\ref{pro:extreme_fillings} that a skew shape $A$ has an
LR-filling with content $\rows{}{A}$ and that this is its least dominant LR-filling.  
We now give a direct description of this filling, since it will be useful later.
First consider the rightmost box of each non-empty row of $A$, and fill these with the numbers
$1, 2, \ldots$ from top to bottom.  Then apply this procedure to the skew shape consisting
of the boxes of $A$ that have not yet been filled.  Repeat until every box of $A$ has been filled.
It is readily checked that the resulting filling is an LR-filling with content $\rows{}{A}$.
See Figure~\ref{fig:extreme_fillings}(b) for an example.   
\end{remark}

Now that we have discussed the well-known necessary conditions, we are ready to
describe our new necessary conditions for Schur-positivity.  
These conditions are inspired
by the necessary conditions for skew Schur equality of \cite{RSV07} and we 
begin with the relevant background from \cite{RSV07}.  The central definition
gives a measure of the amount of overlap among the rows of a skew shape and 
among the columns.

\begin{definition}
Let $A$ be a skew shape with $r$ rows.  For $i = 1, \ldots, r-k+1$, define $\overlap{k}{i}$
to be the number of columns occupied in common by rows $i, i+1, \ldots, i+k-1$.  
Then $\rows{k}{A}$ is defined to be the weakly decreasing rearrangement of 
$(\overlap{k}{1}, \overlap{k}{2}, \ldots, \overlap{k}{r-k+1})$.  
Similarly, we define $\cols{k}{A}$ by looking at the overlap among the columns of $A$.  
\end{definition}

In particular, note that $\rows{1}{A} = \rows{}{A}$ and $\cols{1}{A} = \cols{}{A}$.  

\begin{example}
Let $A=553111/31$ as in Figure~\ref{fig:extreme_fillings}.  We have that 
$\rows{1}{A} = 432111$, $\rows{2}{A}=22111$, $\rows{3}{A}=11$, $\rows{4}{A}=1$ and
$\rows{i}{A} = 0$ otherwise.  Also
$\cols{1}{A} = 42222$, $\cols{2}{A}=2211$, $\cols{3}{A}=111$, $\cols{4}{A}=1$ and
$\cols{i}{A} = 0$ otherwise.
\end{example}

The following result is taken directly from \cite{RSV07}.

\begin{proposition}
Given a skew shape $A$, consider the doubly-indexed array 
\[
(\rects{k}{l}{A})_{k,l \geq 1}
\]
where $\rects{k}{l}{A}$ is defined to be the number of $k \times l$ 
rectangular subdiagrams contained inside
$A$.  Then we have
\begin{eqnarray}
\rects{k}{l}{A} & = & \sum_{l' \geq l} \left( \rows{k}{A} ^t \right)_{l'} \label{equ:rects} \\
& = & \sum_{k' \geq k} \left( \cols{l}{A} ^t \right)_{k'} . \nonumber
\end{eqnarray}
Consequently, any one of the three forms of data
\[
(\rows{k}{A})_{k \geq 1},\ \  (\cols{l}{A})_{l \geq 1}, \ \  (\rects{k}{l}{A})_{k, l \geq 1}
\]
on $A$ determines the other two uniquely.  
\end{proposition}

Note that \eqref{equ:rects} tells us that $\rects{k}{l}{A}$ is the number of boxes weakly to the 
right of column $l$ in the partition $\rows{k}{A}$.  

Not only do $(\rows{k}{A})_{k \geq 1}$, $(\cols{l}{A})_{l \geq 1}$ and $(\rects{k}{l}{A})_{k, l \geq 1}$
determine each other, but their inequalities are related in the following sense.

\begin{proposition}\label{pro:3measures}
Let $A$ and $B$ be skew shapes.  Then $\rows{k}{A} \domleq \rows{k}{B}$ if and only if 
$\rects{k}{l}{A} \leq \rects{k}{l}{B}$ for all $l$.  Similarly, $\cols{l}{A} \domleq \cols{l}{B}$ if and only if 
$\rects{k}{l}{A} \leq \rects{k}{l}{B}$ for all $k$.  Consequently, the following are equivalent:
\begin{itemize}
\item $\rows{k}{A} \domleq \rows{k}{B}$ for all $k$;
\item $\cols{l}{A} \domleq \cols{l}{B}$ for all $l$;
\item $\rects{k}{l}{A} \leq \rects{k}{l}{B}$ for all $k, l$.  
\end{itemize}
\end{proposition}

\begin{proof}
Suppose that $\rows{k}{A} \domleq \rows{k}{B}$ for some fixed $k$.  For any fixed $l$, we wish
to show that $\rects{k}{l}{A} \leq \rects{k}{l}{B}$; i.e.\ the number of elements of the partition
$\rows{k}{A}$ weakly to the right of column $l$ is less than or equal to the number of elements
of $\rows{k}{B}$ weakly to the right of column $l$.  
Suppose column $l$ in $\rows{k}{A}$ and 
$\rows{k}{B}$ has length $a$ and $b$ respectively.  
If $a \leq b$ then we have 
\begin{eqnarray*}
\rects{k}{l}{A}  & =  & \sum_{i=1}^{a} \left((\rows{k}{A})_i -l+1\right) \\
& = & \left(\sum_{i=1}^{a} (\rows{k}{A})_i\right) -a(l-1) \\
& \leq & \left(\sum_{i=1}^{a} (\rows{k}{B})_i\right) -a(l-1)\\
& = & \sum_{i=1}^{a} \left((\rows{k}{B})_i -l+1\right) \\
& \leq & \sum_{i=1}^{b} \left((\rows{k}{B})_i -l+1\right) \\
& = & \rects{k}{l}{B}.
\end{eqnarray*}
If $a > b$ then
\begin{eqnarray*}
\rects{k}{l}{A}  & =  & \sum_{i=1}^{a} \left((\rows{k}{A})_i -l+1\right) \\
& = & \left(\sum_{i=1}^{a} (\rows{k}{A})_i\right) -a(l-1) \\
& \leq & \left(\sum_{i=1}^{a} (\rows{k}{B})_i\right) -a(l-1)\\
& = & \left(\sum_{i=1}^{b} (\rows{k}{B})_i\right) + \left(\sum_{i=b+1}^{a} (\rows{k}{B})_i\right)  -a(l-1)\\
& \leq & \left(\sum_{i=1}^{b} (\rows{k}{B})_i\right) + (a-b)(l-1) - a(l-1) \\
& = & \left(\sum_{i=1}^{b} (\rows{k}{B})_i\right) -b(l-1) \\
& = & \rects{k}{l}{B}.
\end{eqnarray*}

Now suppose $\rects{k}{l}{A} \leq \rects{k}{l}{B}$ for all $l$.  For any fixed $j$, we wish to show that
\[
\sum_{i=1}^{j} (\rows{k}{A})_i \leq \sum_{i=1}^{j} (\rows{k}{B})_i\ .
\]
Suppose row $j$ in $\rows{k}{A}$ and 
$\rows{k}{B}$ has length $a$ and $b$ respectively.  
If $a \leq b$ then, referring to Figure~\ref{fig:3measures}(a), we see that
\begin{eqnarray}
\sum_{i=1}^{j} (\rows{k}{A})_i & =  & ja + \rects{k}{a+1}{A} \label{equ:topleft} \\
& \leq & jb + \rects{k}{b+1}{A} \label{equ:bottomleft}\\
& \leq & jb + \rects{k}{b+1}{B} \nonumber \\
& = & \sum_{i=1}^{j} (\rows{k}{B})_i\ . \nonumber
\end{eqnarray}

\begin{figure}[htpb]
\[
\scalebox{.6}{\includegraphics{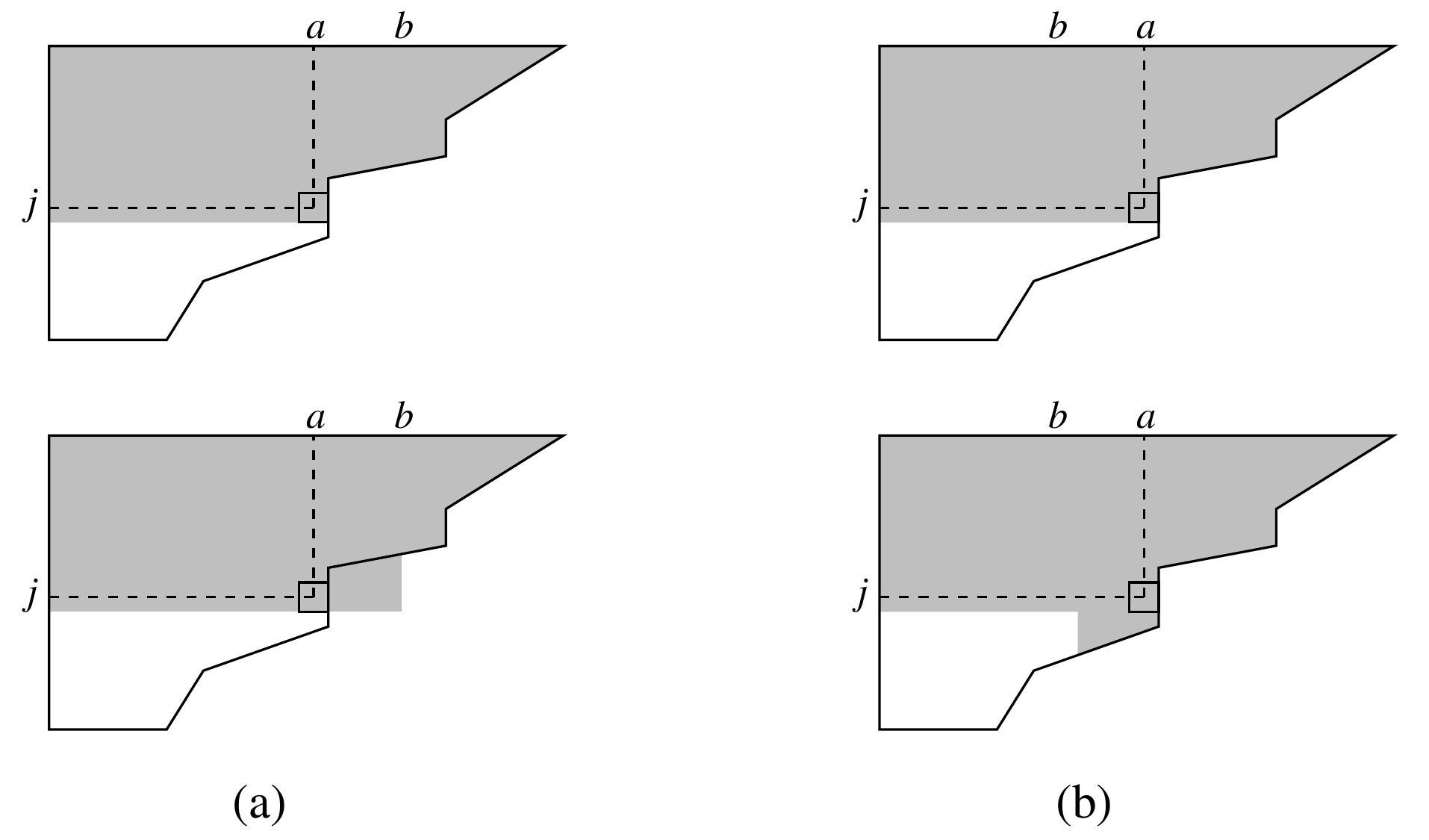}}
\]
\caption{Demonstration of inequalities in the proof of Proposition~\ref{pro:3measures}.  All diagrams represent $\rows{k}{A}$.  The shaded regions in the diagrams on the left represent 
\eqref{equ:topleft} and \eqref{equ:bottomleft}.  Likewise, the shaded regions on the right represent
\eqref{equ:topright} and \eqref{equ:bottomright}.}
\label{fig:3measures}
\end{figure}

If $a > b$ then, referring to Figure~\ref{fig:3measures}(b), we see that
\begin{eqnarray}
\sum_{i=1}^{j} (\rows{k}{A})_i & = & ja + \rects{k}{a+1}{A} \nonumber \\
& = & jb + j(a-b) + \rects{k}{a+1}{A} \label{equ:topright}\\
& \leq & jb + \rects{k}{b+1}{A} \label{equ:bottomright} \\
& \leq & jb + \rects{k}{b+1}{B} \nonumber \\
& = & \sum_{i=1}^{j} (\rows{k}{B})_i\ . \nonumber
\end{eqnarray}

We conclude that  $\rows{k}{A} \domleq \rows{k}{B}$ if and only if 
$\rects{k}{l}{A} \leq \rects{k}{l}{B}$ for all $l$.
The proof that $\cols{l}{A} \domleq \cols{l}{B}$ if and only if 
$\rects{k}{l}{A} \leq \rects{k}{l}{B}$ for all $k$ is similar.  It is then an easy consequence that
the three sets of inequalities are equivalent.  
\end{proof}

For any skew shape $A$, we let $\trim{}{A}$ denote the skew shape obtained from $A$ by deleting
the top element of every non-empty column of $A$.  
We will consider $\trim{}{A}$ to be a function on skew shapes, so that 
$\trim{k}{A} = \trim{}{\trim{k-1}{A}}$ and $\trim{1}{A}$ is simply $\trim{}{A}$.
 
\begin{lemma}\label{lem:trim}  Let $A$ be any skew shape and let $k$ be an integer with
$k \geq 2$.  Then 
\[
\rows{k-1}{\trim{}{A}} = \rows{k}{A}.
\]
\end{lemma}

\begin{proof}
Note that if the $i$th row of $A$ has $c$ columns in common with the $(i+k-1)$st row of $A$, then
the $(i+1)$st row of $\trim{}{A}$ has exactly $c$ columns in common with the $(i+k-1)$st row of
$\trim{}{A}$.  The result follows.
\end{proof}

We are now in a position to state and prove the central part of our main result.

\begin{theorem}\label{thm:main}
Let $A$ and $B$ be skew shapes.
If $\supp{A} \supseteq \supp{B}$, 
then 
\[
\rows{k}{A} \domleq \rows{k}{B} \mbox{\ for all $k$}.
\]
\end{theorem}

\begin{proof}
We consider a particular LR-filling of $B$.  Roughly speaking, we will fill $B$ with the numbers 
$1,2,\ldots,k-1$ in the most dominant way possible, and then fill the boxes that remain with
the numbers $k, k+1, \ldots$ in the least dominant way possible.  
More precisely, fill the $i$th
highest box of each column of $B$, when such a box exists, with $i$, for $i=1,2,\ldots, k-1$.  
Thus there will be $(\cols{}{B}^t)_i$ boxes of $B$ filled with $i$, for $i=1,2,\ldots,k-1$.
See Figure~\ref{fig:main_proof} for an example of the entire proof in action in the case $k=3$.
\begin{figure}[htpb]
\[
\scalebox{.7}{\includegraphics{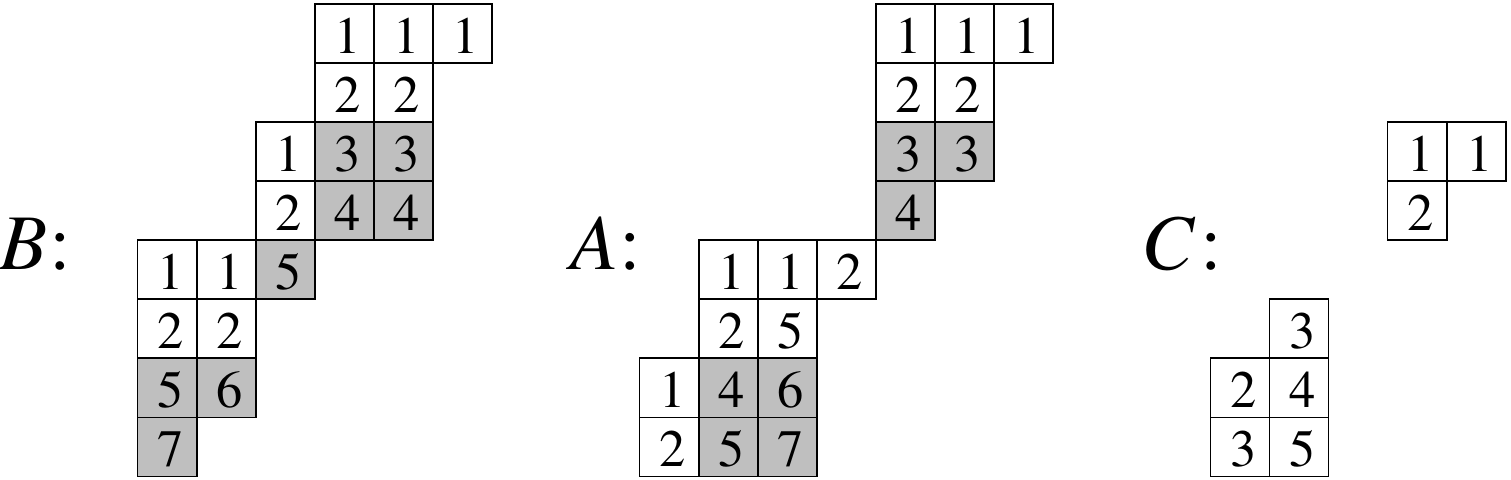}}
\]
\caption{Demonstration of the proof of Theorem~\ref{thm:main} with $k=3$; $\trim{k-1}{B}$ and $\trim{k-1}{A}$ are shaded}
\label{fig:main_proof}
\end{figure}
We see that the boxes of $B$ that remain empty form the skew shape $\trim{k-1}{B}$.  
We fill these boxes with
the numbers $k, k+1, \ldots$ in the least dominant way possible, as described in 
Remark~\ref{rem:extreme_fillings}.  Thus there will be $\rows{}{\trim{k-1}{B}}_{i-k+1}$ boxes of $B$ filled
with $i$, for $i=k, k+1, \ldots$.  It is straightforward to check that the overall result is an LR-filling
of $B$.  The key observation is that, by repeated applications of 
Lemma~\ref{lem:trim}, $\rows{}{\trim{k-1}{B}} = \rows{k}{B}$.

Since $\supp{B} \subseteq \supp{A}$, there must be an LR-filling of $A$ with content
\[
\left((\cols{}{B}^t)_1, 
\ldots,  (\cols{}{B}^t)_{k-1}, \rows{}{\trim{k-1}{B}}_{1}, \rows{}{\trim{k-1}{B}}_{2}, \ldots \right).
\]
In our running example, Figure~\ref{fig:main_proof} shows one possibility.
Removing the boxes of $A$ filled with $1,2, \ldots, k-1$ in this filling results in a skew shape
$C$ that is filled with the numbers $k, k+1, \ldots$.  We see that subtracting $k-1$ from the 
entries of the boxes of $C$ results in an LR-filling of $C$ of content $\rows{}{\trim{k-1}{B}}$.  
This is the filling of $C$ shown in Figure~\ref{fig:main_proof}.
By Proposition~\ref{pro:extreme_fillings}, we deduce that $\rows{}{C} \domleq \rows{}{\trim{k-1}{B}}$.

Now consider $\trim{k-1}{A}$.  
As with $\trim{k-1}{B}$, by repeated
applications of Lemma~\ref{lem:trim}, $\rows{}{\trim{k-1}{A}} = \rows{k}{A}$.  Also note that 
in any SSYT of shape
$A$, the numbers $1,2,\ldots,k-1$ can only appear in the top $k-1$ boxes of some column of $A$.
Therefore, $\trim{k-1}{A} \subseteq C$, by definition of $C$.  Applying Lemma~\ref{lem:dom_containment},
we deduce that $\rows{}{\trim{k-1}{A}} \domleq \rows{}{C}$, and so 
$\rows{}{\trim{k-1}{A}} \domleq \rows{}{\trim{k-1}{B}}$.  This is exactly the desired inequality:
$\rows{k}{A} \domleq \rows{k}{B}$.
\end{proof}

As a special case of Theorem~\ref{thm:main}, we achieve our main goal of obtaining 
necessary conditions for Schur-positivity.

\begin{corollary}\label{cor:Schurpositivity}
Let $A$ and $B$ be skew shapes.
If $s_A - s_B$ is Schur-positive, 
then 
\[
\rows{k}{A} \domleq \rows{k}{B} \mbox{\ for all $k$}.
\]
\end{corollary}

Combining Theorem~\ref{thm:main} and Corollary~\ref{cor:Schurpositivity} 
with Proposition~\ref{pro:3measures}, we actually get three equivalent sets of necessary
conditions for support containment or Schur-positivity.

\begin{corollary}\label{cor:combine}
Let $A$ and $B$ be skew shapes.
If $s_A - s_B$ is Schur-positive, 
or if $A$ and $B$ satisfy the weaker condition that $\supp{A} \supseteq \supp{B}$, then
the following three equivalent conditions are true:
\begin{itemize}
\item $\rows{k}{A} \domleq \rows{k}{B}$ for all $k$;
\item $\cols{l}{A} \domleq \cols{l}{B}$ for all $l$;
\item $\rects{k}{l}{A} \leq \rects{k}{l}{B}$ for all $k, l$.  
\end{itemize}
\end{corollary}

\begin{example}
Let 
\[
\begin{picture}(4,4)(1,0)
\put(1,4){\line(1,0){3}}
\put(0,3){\line(1,0){4}}
\put(0,2){\line(1,0){2}}
\put(0,1){\line(1,0){1}}
\put(0,0){\line(1,0){1}}
\put(0,3){\line(0,-1){3}}
\put(1,4){\line(0,-1){4}}
\put(2,4){\line(0,-1){2}}
\multiput(3,4)(1,0){2}{\line(0,-1){1}}
\put(-2,2){$A=$}
\end{picture}
\begin{picture}(4,4)(-4.5,0)
\put(-4.5,2){and\quad$B=$}
\put(3,4){\line(1,0){1}}
\put(1,3){\line(1,0){3}}
\put(1,2){\line(1,0){3}}
\put(0,1){\line(1,0){3}}
\put(0,0){\line(1,0){1}}
\put(4,4){\line(0,-1){2}}
\put(3,4){\line(0,-1){3}}
\put(2,3){\line(0,-1){2}}
\put(1,3){\line(0,-1){3}}
\put(0,1){\line(0,-1){1}}
\put(4.5,2){.}
\end{picture}
\]
We see that
$\rows{2}{A} = 111$ and $\rows{2}{B} = 21$.  Thus we know that $s_B - s_A$ is 
not Schur-positive.  On the other hand, $\rows{3}{A} = 1$ while $\rows{3}{B} = 0$, 
implying that $s_A - s_B$ is not Schur-positive.  Moreover, we can conclude that
$\supp{A}$ and $\supp{B}$ are incomparable under containment order.

It is certainly not the case that $\rows{k}{A} \domleq \rows{k}{B}$ for all $k$ implies that
$\supp{A} \supseteq \supp{B}$.  To see the smallest example, let $A = 311/1$ and 
$B=22$, and note that $22 \in \supp{B} \setminus \supp{A}$.
\end{example}

\section{Special Cases}\label{sec:specialcases}

\subsection{Strengthening of skew Schur equality result}

The central result of the penultimate section of \cite{RSV07}, namely Corollary~8.11, states
that if $s_A = s_B$ then $\rows{k}{A} = \rows{k}{B}$ for all $k$.  Since $\domleq$ is a partial
order, we can use Theorem~\ref{thm:main} to 
replace the hypothesis $s_A = s_B$ with the weaker hypothesis 
$\supp{A} = \supp{B}$ as follows.

\begin{corollary}
Let $A$ and $B$ be skew shapes.
If $\supp{A} = \supp{B}$, then 
\[
\rows{k}{A} = \rows{k}{B} \mbox{\ for all $k$}.
\] 
\end{corollary}

\begin{example}
We note that pairs of skew shapes $(A, B)$ with $\supp{A} = \supp{B}$ but $s_A \neq s_B$
are quite common.  When $|A| = |B| = 5$, there are already several examples, including
$A = 3311/21$ and $B=3321/211$.  We see that 
\[
s_A = s_{32} + s_{2111} + s_{221} + s_{311} \mbox{\ \ \ and\ \ \ } s_B = s_{32} + s_{2111} + 2s_{221} + s_{311}.
\]
\end{example}

\subsection{Products of Schur functions}\label{sub:products}

\begin{figure}[htpb]
\[
\scalebox{.5}{\includegraphics{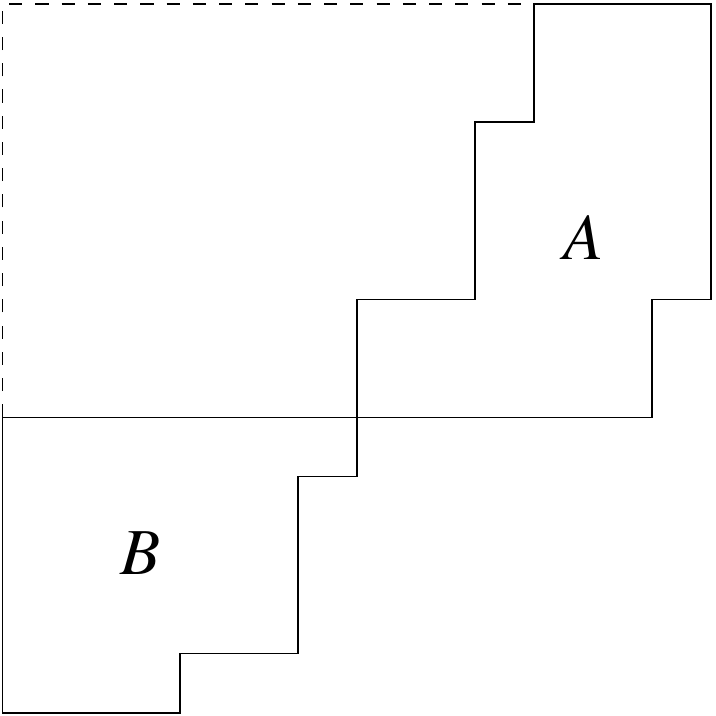}}
\]
\caption{Positioning two skew shapes $A$ and $B$ as shown results in another skew shape}
\label{fig:product}
\end{figure}

Any product $s_A s_B$ of skew Schur functions $s_A$ and $s_B$ 
is again a skew Schur function, as made evident by Figure~\ref{fig:product} and the definition
of skew Schur functions.  We denote the skew shape of Figure~\ref{fig:product} by $A \ast B$.
In particular, for partitions $\gamma$ and $\delta$, the product $s_\gamma s_\delta$ of 
Schur functions is a skew Schur function.  The Schur-positivity of the expression 
\begin{equation}\label{equ:fflp}
s_\alpha s_\beta - s_\gamma s_\delta
\end{equation}
for partitions $\alpha$, $\beta$ has been studied, for example, in \cite{BBR06,FFLP05,LPP07}.
It is natural to ask what Theorem~\ref{thm:main} tells us about expressions of the form 
\eqref{equ:fflp}.  We begin with the following lemma that is simple to check.  

\begin{lemma}\label{lem:straight}
For any straight shape $\alpha = (\alpha_1, \alpha_2, \ldots, \alpha_l)$, we have
\[
\rows{k}{\alpha} = (\alpha_k, \alpha_{k+1}, \ldots, \alpha_l)
\]
 for $k=1,2,\ldots,l$.
\end{lemma}

For partitions $\alpha$ and $\beta$, let $\alpha \union \beta$ denote the partition whose multiset
of parts equals the union of the multisets of parts of $\alpha$ and $\beta$.
The following result is essentially Theorem~\ref{thm:main}
specialized to skew shapes of the form $A \ast B$.

\begin{corollary}\label{cor:product}
If partitions $\alpha$, $\beta$, $\gamma$ and $\delta$ satisfy the condition that
$s_\alpha s_\beta - s_\gamma s_\delta$ is Schur-positive,
or satisfy the weaker condition that 
$\supp{s_\alpha s_\beta} \supseteq \supp{s_\gamma s_\delta}$,
then 
\[
(\alpha_k, \alpha_{k+1},\ldots, \alpha_l) \union
(\beta_k, \beta_{k+1},\ldots, \beta_l)  \domleq
(\gamma_k, \gamma_{k+1},\ldots, \gamma_l) \union
(\delta_k, \delta_{k+1},\ldots, \delta_l)
\]
for all $k$, and for all 
$l \geq \max\{\ell(\alpha), \ell(\beta)\}$. 
\end{corollary}

\begin{proof}
We know that $s_\alpha s_\beta = s_{\alpha \ast \beta}$.  It is also clear that 
$\rows{k}{\alpha \ast \beta} = \rows{k}{\alpha} \union \rows{k}{\beta}$.  The result now
follows from Theorem~\ref{thm:main} and Lemma~\ref{lem:straight}.
\end{proof}

In words, we might say that all the ``tails'' from $\alpha$ and $\beta$ are dominated by those from
$\gamma$ and $\delta$.


%
%



\bibliography{necconditions}
\bibliographystyle{plain}

\end{document}